\documentclass[a4paper]{siamltex}

\usepackage[english]{babel}
\usepackage{amssymb,amsfonts,amsmath,latexsym}
\usepackage[dvips]{graphicx}
\usepackage{algorithm}
\usepackage{algorithmic}
\usepackage{color}
\usepackage{url}
\usepackage{hyperref}

\usepackage[margin=1.3in]{geometry}

\newcommand{\bm}[1]{\mathbf{#1}}
\newcommand{\R}{\mathbb{R}}

\newcommand{\cS}{{\mathcal{S}}}

\newcommand{\ellb}{\boldsymbol{\ell}}

\begin{document}

\title{Ascertaining the ideality of photometric stereo datasets under unknown
lighting
}

\author{Elisa Crabu\thanks{Department of Mathematics and Computer Science, 
University of Cagliari, via Ospedale 72, 09124 Cagliari, Italy, 
\texttt{elisacrabu96@outlook.it, rodriguez@unica.it}}
\and
Federica Pes\thanks{Department of Chemistry  and Industrial Chemistry, 
University of Pisa, Via G. Moruzzi, 13, 56124 Pisa, Italy, 
\texttt{federica.pes@dcci.unipi.it}}
\and
Giuseppe Rodriguez\footnotemark[1] 
\and
Giuseppa Tanda\thanks{CeSim -- Centro Studi ``Identit\`a e Memoria'', 07100 
Sassari, Italy, \texttt{giuseppa.tanda@gmail.com}}
}

\maketitle

\medskip

\begin{abstract}
The standard photometric stereo model makes several assumptions that
are rarely verified in experimental datasets. In particular, the observed
object should behave as a Lambertian reflector and the light sources should be
positioned at an infinite distance from it, along a known direction. Even when
Lambert's law is approximately fulfilled, an accurate assessment of the
relative position between the light source and the target is often unavailable
in real situations. The Hayakawa procedure is a computational method for
estimating such information directly from the data images. It occasionally
breaks down when some of the available images {excessively deviate from
ideality.
This is generally due to observing a non Lambertian surface, or illuminating it
from a close distance, or both.
Indeed, in narrow shooting scenarios, typical, e.g., of archaeological
excavation sites, it is impossible to position a flashlight at a sufficient
distance from the observed surface.} It is then necessary to understand if a
given dataset is reliable and which images should be selected to better
reconstruct the target. In this paper, we propose some algorithms to perform
this task and explore their effectiveness.
\end{abstract}



\section{Introduction}

Photometric stereo (which we will denote by PS) is a classic computer vision
approach for reconstructing the shape of a three-dimensional
object~\cite{Woodham1979,Woodham1980}.
It is considered a \emph{shape from shading} (SfS) technique, as it takes as
input images that embed shape and color information of the observed object. 
However, while the original SfS problem only considers a single two-dimensional
image of the object~\cite{dfs2008,rpjm1999} and it is known to be ill-posed, PS
stands on the use of a set of images acquired from a fixed point of view under
varying lighting conditions.

If the lighting positions are known,
Kozera~\cite{kozera1991} proved that, under suitable assumptions, the PS
problem is well-posed when at least two pictures are available.
Kozera's approach consists of modeling the problem by a system of first-order
Hamilton-Jacobi partial differential equations (PDEs), for which various
numerical methods have been proposed; see, e.g.,~\cite{mecca2013}.
The solution process that we will consider is slightly different and operates
in two steps.
First, the normal vectors at each discretization point of a rectangular
domain, namely, the digital picture, are determined by solving an algebraic
matrix equation.
Then, after numerically differentiating the normal vector field, the solution
of a Poisson PDE leads to the approximation of the object surface.
This approach has been used in~\cite{dmrtv15} to reconstruct rock art carvings
found in \emph{Domus de Janas} (fairy houses), a specific typology of Neolithic
tombs found in Sardinia, Italy~\cite{vmdrt14}.
While this computational scheme requires at least three images with different
lighting conditions, the availability of a larger dataset leads to
least-squares solution method, which allows for a more efficient treatment of
measurement errors.

Though in the research community the position of the light sources is generally
assumed to be known~\cite{barsky2003,khanian2018,Radow2018}, real-world
applications of PS often deal with datasets acquired under unknown light
conditions~\cite{basri2007,chen2006}.
In~\cite{hayakawa1994}, Hayakawa has shown that the lighting directions can be
identified directly from the data when at least six images with different
illumination are available {(see~\cite{hayakawa1994} and
Section~\ref{sec:hayakawa} for a proof)}, yielding the possibility to apply PS
to field measurement scenarios, like in the case of archaeological
excavations~\cite{cdfrv21}.

\begin{figure}[htb]
\centerline{\includegraphics[width=\linewidth]{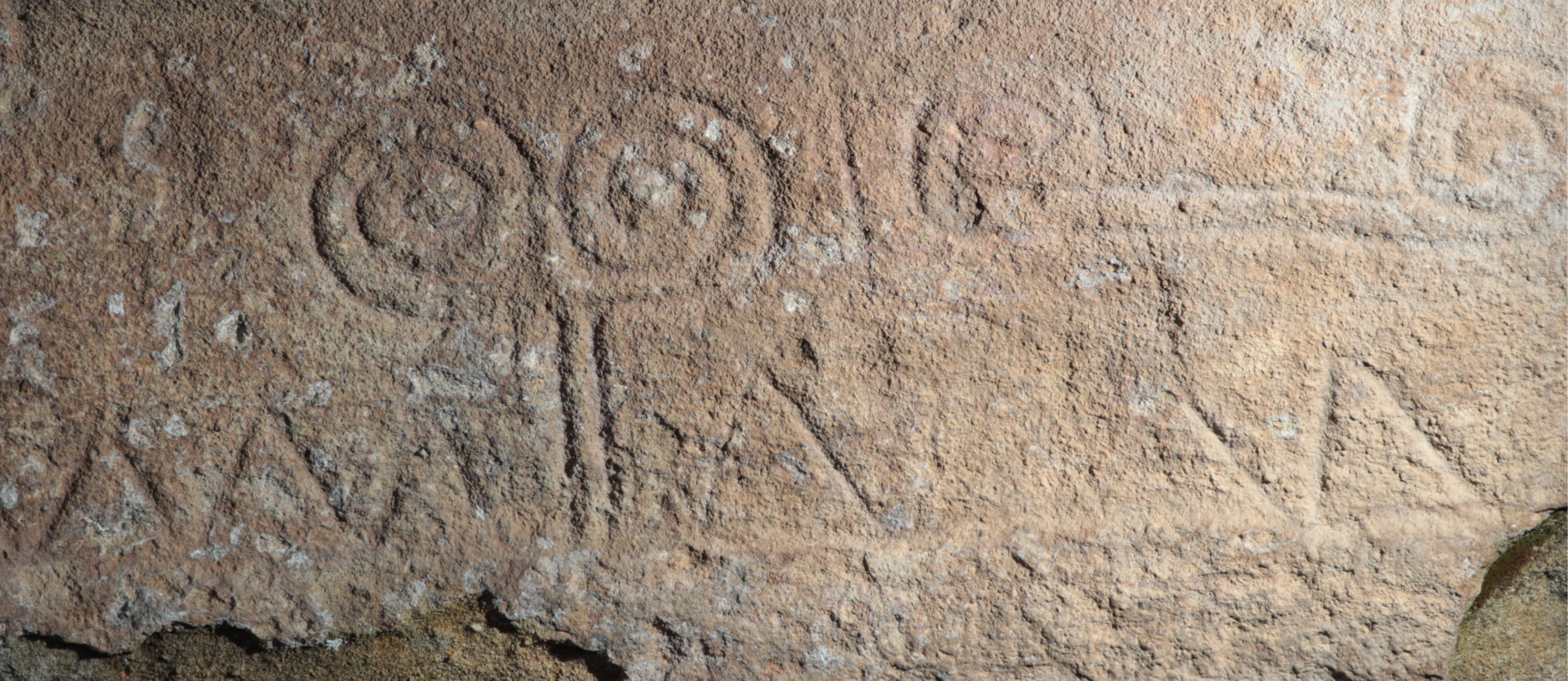}}
\caption{Engraving found in the \emph{Domus de Janas} of Corongiu, in Pimentel
(Sardinia, Italy).}
\label{fig:corongiu}
\end{figure}

\begin{figure}[htb]
\centerline{\includegraphics[width=.8\linewidth]{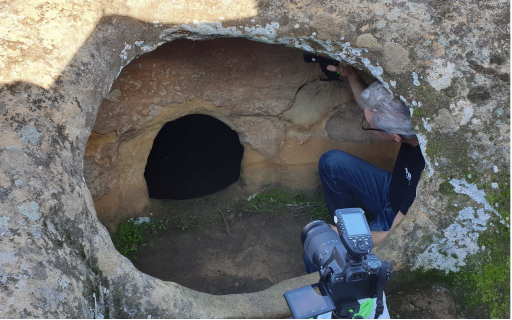}}
\caption{Lobby of the Neolithic tomb located in the \emph{Domus de Janas} of
Corongiu (Sardinia, Italy). The engraving is over the entrance to the internal
chamber.}
\label{fig:pimentel}
\end{figure}

For example, the bas-relief displayed in Figure~\ref{fig:corongiu} is found in
the \emph{Domus de Janas} of Corongiu, in Pimentel (South Sardinia, Italy),
dated 4th millennium BC.
Figure~\ref{fig:pimentel} shows the small lobby where the engraving is located.
While the camera can be placed sufficiently far from the rock surface,
it is impossible to illuminate it with a flashlight at a distance
suitable for proper lighting. Moreover, measuring with precision the 
position of the light source relative to the engraving is very difficult.
This situation is common in many archaeological sites.

PS application implies that certain assumptions are met in the acquisition
process.
In general, a perspective camera projection model may be used, with parameters
estimated by camera calibration methods; see \cite{YaBr2020} for an
introduction to the use of perspective and~\cite{khanian2018,mtw2014} for
applications.
In this work, we consider an orthographic camera model, which simplifies the
presentation and is suitable for situations where the camera can be positioned
at a relatively far distance from the observed object.
The mathematical formulation of PS assumes that the object surface is
Lambertian, meaning that Lambert's cosine law can be employed to
describe the object reflectance.
This condition implies that the surface is matte and free from specular
reflections, which is seldom verified in practice; see~\cite{tmdd2016,tf2016}
and~\cite{rrmyb23}, where the Oren--Nayar model has been used to preprocess a
dataset originated by a non-Lambertian object.

Many studies have been devoted to improving the application of PS to
non-Lambertian conditions.
In~\cite{wu2011robust}, the authors formulate the problem as a low-rank matrix
factorization subject to sparse errors, due to the presence of shadows and
specular reflections in data images.
A statistical approach based on a hierarchical Bayesian approximation was used
in~\cite{ikehata2012robust} to simultaneously estimate the surface normal
vectors and the experimental errors, by solving a constrained sparse regression
problem.
Learning-based methods were also employed to solve PS; see,
e.g.,~\cite{logothetis2021px,ju2022normattention}. 
To manage near-light problems, models based on neural
networks were presented in~\cite{logothetis2023cnn,lichy2022fast}, where
they were stated under non-Lambertian conditions.

Our approach is slightly different from the ones above.
Here, we assume that the light source location is unknown, and point our
attention to the lack of ideality in a dataset due to either the necessity of
illuminating the observed surface with close lights, or to the fact that the
surface is only approximately Lambertian.
Indeed, in archaeological applications, which motivated this work, the
assumption which requires the light sources to be located at a large distance
(theoretically, at infinite distance) from the observed surface is an
unacceptable constraint, as carvings and bas-reliefs are often located in
narrow caves or excavations.
Furthermore, the carvings are engraved in rock, that may not be an ideal Lambertian
surface.
To explore these aspects of the model, we assume that orthographic projection
conditions are met and the position of light sources is unknown.
Then, we determine the light position by the Hayakawa procedure and
investigate the reason for its breakdown in the presence of nonideal data,
i.e., images that excessively deviate from the above assumption.
This also leads us to introduce a new nonlinear approach for identifying
the light position, alternative to Hayakawa's one.
Finally, we study the problem of ascertaining the ideality of a dataset. This
is accomplished by defining some measures of ideality, which can assist a researcher in
selecting the best subset of images from a redundant dataset.
To demonstrate the performance of the new approaches, we apply them to both a
synthetic dataset, where the position of the lights can be decided at will,
and to an experimental one, which is not guaranteed to represent a Lambertian
reflector.

The plan of the paper is the following. Section~\ref{sec:hayakawa} contains a
review of the numerical procedures to solve the PS problem under both known and
unknown lighting.
In Section~\ref{sec:nonlin}, we describe a new nonlinear approach to determine
the light position directly from the available data.
The ideality of a dataset is discussed in Section~\ref{sec:nonideal}, where two
indicators are introduced for detecting the presence of unideal images, and an
algorithm is presented to extract a subset of pictures that better suit the
assumptions of the PS model.
The numerical experiments that investigate the performance of the algorithm are
outlined in Section~\ref{sec:expresults} and discussed in
Section~\ref{sec:discussion}.

\section{A review of the Hayakawa procedure for determining the light
position}\label{sec:hayakawa}

In this section, we briefly review the mathematical model implementing
Lambert's law that is usually adopted in photometric stereo under the
assumption of orthographic projection and light source at an infinite distance
from the target; see~\cite{cdfrv21} for more details.
We also outline the Hayakawa procedure~\cite{hayakawa1994} for estimating the
lighting directions from the data.

Let us associate an orthonormal reference system to $\mathbb{R}^3$ and
assume the observed object is located at its origin.
The optical axis of the camera coincides with the $z$-axis and the point of
view is at infinite distance, to ensure orthographic projection.
The camera has a resolution of $(r+2)\times(s+2)$ pixels, 
and is associated with the rectangular domain
$\Omega = \left[ -A/2,A/2 \right] \times \left[ -B/2,B/2 \right]$, 
where $A$ is the horizontal size of the rectangle and $B=(s+1)h$, with
$h=A/(r+1)$. 
The domain discretization consists of a grid of points with coordinates
$(x_i,y_j)$, for $i=0,\dots,r+1$ and $j=0,\dots,s+1$, given by
\begin{equation}\label{discretization}
x_i = -\frac{A}{2} + ih, \qquad
y_j = -\frac{B}{2} + jh. 
\end{equation}

The surface of the object is represented by a bivariate explicit function
$z=u(x,y)$, $(x,y)\in \Omega$, whose partial derivatives are denoted by $u_x$
and $u_y$.
Then, the gradient of $u$ and the (normalized) normal vector to each point of
the surface are defined as
\begin{equation*}
\nabla  u(x,y) = \begin{bmatrix}  u_x \\ u_y \end{bmatrix}
\qquad \text{and} \qquad
\bm{n}(x,y) = \frac{[-u_x,-u_y,1]^T}{\sqrt{1 + \| \nabla u\|^2}},
\end{equation*}
respectively, where $\|\cdot\|$ denotes the 2-norm.
All functions are discretized on the grid, and their values are stacked in a
vector in lexicographic order, that is, according to the rule $k=(i-1)s+j$,
where $k=1,\dots,p$, and $p$ denotes the number of pixels in the image.
We will write indifferently $u(x_i,y_j)$, $u_{i,j}$, or $u_k$, and similarly
for $u_x(x_i,y_j)$, $u_y(x_i,y_j)$, and $\bm{n}(x_i,y_j)$.

The images in the dataset are stacked in the vectors $\bm{m}_t\in\mathbb{R}^p$,
$t=1,\ldots,q$. In each of the pictures, the target is illuminated by a light
source located at a different direction.
We denote the vector that points from the object to the light source by
$\ellb_t=[\ell_{1t},\ell_{2t},\ell_{3t}]^T$, with $t=1,\dots,q$.
Its 2-norm is proportional to the light intensity. 
To simplify notation, in the following we assume the light vectors to be
normalized, that is, $\|\ellb_t\|=1$.

Assuming that the observed surface is a Lambertian reflector, we state
Lambert's cosine law in the form
\begin{equation}\label{eqLambert}
\rho(x,y) \langle \bm{n}(x,y), \ellb_t \rangle = I_t(x,y), \qquad t=1,\dots,q,
\end{equation}
where $\langle\cdot,\cdot\rangle$ is the standard inner product in $\R^3$ and
the albedo $\rho(x,y)$ accounts for the partial light absorption of the
surface, due to its color or material.
After discretization, formula \eqref{eqLambert} reads
\begin{equation}\label{discLambert}
    \rho_k \bm{n}_k^T \ellb_t = m_{kt}, \qquad k=1, \dots, p, \  t=1, \dots, q,
\end{equation}
where $m_{kt}$ represents the $k$th component of the image vector $\bm{m}_t$,
that is, the radiation $I_t(x_i,y_j)$ reflected by a neighborhood of the $k$th
pixel in the image number $t$.

Equations \eqref{discLambert} can be represented in the form of the matrix
equation
\begin{equation}\label{problem}
D N^T L = M,
\end{equation}
after assembling data and unknowns in the matrices
\begin{equation*}
 \begin{aligned}
    D & = \diag(\rho_1, \dots, \rho_p) \in \mathbb{R}^{p \times p},
    &\quad
    L &= \left[ \ellb_1, \dots, \ellb_q\right] \in\mathbb{R}^{3 \times q}, \\
    N &=\left[ \bm{n}_1, \dots, \bm{n}_p \right] \in \mathbb{R}^{3 \times p},
    &\quad
    M &=\left[ \bm{m}_1, \dots, \bm{m}_q \right] \in \mathbb{R}^{p \times q}.
 \end{aligned}
\end{equation*}

If the matrix $L$, whose columns are the vectors aiming at the lights, is
known, one usually sets $\tilde{N}=ND$ in~\eqref{problem}, and computes
$\tilde{N}^T=ML^\dagger$, where $L^\dagger$ denotes the Moore--Penrose
pseudoinverse of $L$~\cite{bjo96}. This produces a unique solution if $q\geq 3$
and the matrix $L$ is full-rank. Once $\tilde{N}$ is available, by normalizing
its columns one easily obtains $N$ and $D$ from which, by integrating the
normal vectors (see~\cite{qda18}), the approximation of the surface
representation $u(x,y)$ at the discretization grid can be computed.

When the light positions are not known, the Hayakawa
procedure~\cite{hayakawa1994} can be used to estimate this information from the
data,  but this requires the availability of at least six images with different
lighting conditions.
Here, we briefly summarize this technique; see~\cite{cdfrv21} for more details.

To obtain an initial rank-3 factorization of the data matrix $M$, we start
computing the singular value decomposition (SVD)~\cite{gvl96} 
\begin{equation}\label{factorization}
M= U \Sigma V^T,
\end{equation}
where $\Sigma = \diag(\sigma_1, \dots, \sigma_q)$ contains
the singular values $\sigma_1\geq\cdots\geq\sigma_q\geq 0$, 
and $U \in \mathbb{R}^{p \times q}$, $V \in
\mathbb{R}^{q \times q}$ are matrices with orthonormal columns $\bm{u}_i$, 
$\bm{v}_i$.
Theoretically, only the first three singular values should be positive, that
is, $\sigma_4=\cdots=\sigma_q=0$. Anyway, because of experimental errors
propagation and lack of ideality in the data, one usually finds $\sigma_i>0$,
$i=4,\ldots,q$.

It is known that the best rank-$k$ approximation of a matrix with respect to
both the 2-norm and the Frobenius norm is produced by truncating the SVD to $k$
terms~\cite{bjo96}.
So, we define $W=[\sigma_1 \bm{u}_1,\sigma_2\bm{u}_2,\sigma_3\bm{u}_3]^T$ and 
$Z=[\bm{v}_1,\bm{v}_2,\bm{v}_3]^T$, to obtain the rank-3 approximation 
$W^TZ$ of $M$.

Following the proof of Theorem 1 from~\cite{cdfrv21}, we consider the
factorization  $W^T Z=M$ as a tentative approximation of $\tilde{N}^TL=M$.
Assuming $\|\ellb_t\|=1$, we seek a matrix $B$ such that 
\begin{equation}\label{Hprob}
\|B\bm{z}_t\|=1, \qquad t=1,\dots,q,
\end{equation}
where $\bm{z}_t\in\R^3$ denote the columns of $Z$.

The Hayakawa procedure consists of writing problem~\eqref{Hprob} in the form
\begin{equation}\label{systemG}
\diag(Z^T G Z)=\mathbf{1},
\end{equation}
where $\mathbf{1}=[1, \dots, 1]^T \in \mathbb{R}^q$ and $G=B^T B$ is a
symmetric positive definite $3 \times 3$ matrix, which depends upon six
parameters, namely, the entries $g_{ij}$ with $i=1,2,3$ and $j=i,\ldots,3$.

Each equation of system \eqref{systemG} reads
\begin{equation*}
\bm{z}_t^T G \bm{z}_t = \sum_{i,j=1}^3 z_{it} z_{jt} g_{ij} = 1.
\end{equation*}
Such equations can be assembled in the linear system $H\bm{g}=\bm{1}$, where 
$H$ is a $q \times 6$ matrix with rows
\begin{equation}\label{matrixH}
\begin{matrix}
[ z_{1t}^2 & z_{2t}^2 & z_{3t}^2 & 2z_{1t}z_{2t} &  2 z_{1t}z_{3t} &
2z_{2t}z_{3t}], \qquad t=1, \dots, q.
\end{matrix}
\end{equation}
The solution $\bm{g}=[g_{11},g_{22},g_{33},g_{12},g_{13},g_{23}]^T$
is unique if the system is overdetermined, that is, {if the matrix $H$ is
full-rank and $q\geq 6$. This shows that a necessary condition for the position
of the light sources to be uniquely determined is that the dataset contains at
least six images.}

The factor $B$ of $G$ is determined up to a unitary transformation, so to
simplify the problem we represent $B$ by its QR factorization and substitute
$B=R$.
It is known~\cite{gvl96} that $R$ can be obtained by the Cholesky factorization
$R^TR$ of the matrix $G$.
This step is particularly important and critical. Indeed, while the matrix $G$
is symmetric by construction, it may be nonpositive definite due to a lack of
ideality in the dataset.
This would prevent the applicability of the Cholesky factorization, causing a
breakdown in the algorithm.

If the previous step is successful, the obtained normal field is usually
rotated, possibly with axes inversions, with respect to the original
orientation of the surface. This would pose some difficulty in the
representation of the surface as an explicit function $z=u(x,y)$ and may lead,
e.g., to a concave reconstruction of a convex surface.
For this reason, the final step is to determine a unitary transformation $Q$
which resolves the uncertainty in the orientation of the normal vectors, the
so-called \emph{bas-relief ambiguity}~\cite{basrelief99}, and rotates the
object to a suitable orientation.
This can be accomplished by an algorithm introduced in~\cite{cdfrv21}, which assumes
the photos are taken by following a particular shooting procedure.

Once $Q$ is determined, the solution of problem~\eqref{problem} is given by
\begin{equation*}
\tilde{N}= QB^{-T}W \quad \text{and} \quad L= QBZ.
\end{equation*}
As already observed, the albedo matrix $D$ and the matrix $N$ containing the
normal vectors are easily obtained by normalizing the columns of $\tilde{N}$.

\section{A nonlinear approach to identify the light position}\label{sec:nonlin}

Here we propose an alternative new method for dealing with
problem~\eqref{Hprob}, based on a nonlinear approach.
By employing, as in the previous section, the QR factorization $B=QR$,
where the matrix $Q$ is orthogonal and $R$ is upper triangular,
we rewrite~\eqref{Hprob} as
\begin{equation}\label{norm2}
\|R\bm{z}_t\|^2 = 1,\qquad t=1,\ldots,q.
\end{equation}
In these equations, the unknowns are 
the nonzero entries $r_{ij}$ of $R$, with $i\leq j$, which we
collect in the vector $\bm{r}=[r_{11}, r_{12}, r_{13}, 
r_{22}, r_{23}, r_{33}]^T$.
Formula~\eqref{norm2} can be seen as a nonlinear equation $F(\bm{r})=\bm{0}$,
for the vector-valued function, with values in $\R^q$, defined by
\begin{equation*}
F(\bm{r}) = [f_1(\bm{r}),\ldots,f_q(\bm{r})]^T,\qquad 
f_t(\bm{r}) = \bm{z}_t^TR^TR\bm{z}_t-1,\qquad t=1,\ldots,q.
\end{equation*}
To determine the solution vector $\bm{r}$, we apply the Gauss--Newton method
to the solution of the nonlinear least-squares problem
\begin{equation}\label{nonlin}
\min_{\bm{r}\in\R^6} \|F(\bm{r})\|^2,
\end{equation}
To ensure the uniqueness of the solution we set $q\geq 6$, the same assumption
of the Hayakawa procedure discussed in the previous section.

We recall that the Gauss--Newton method \cite{bjo96} is an iterative algorithm
that replaces the
nonlinear problem~\eqref{nonlin} by a sequence of linear approximations 
\begin{equation*}
\min_{\bm{s}\in\R^6} \|J(\bm{r}^{(k)})\bm{s} + F(\bm{r}^{(k)})\|,\qquad 
k=0,1,2,\ldots,
\end{equation*}
where $J(\bm{r})$ is the Jacobian matrix of $F(\bm{r})$.
At each iteration, the step $\bm{s}^{(k)}$ is computed by solving the above 
linear least-squares problem. Then, the approximated solution is of the form
\begin{equation*}
\bm{r}^{(k+1)} = \bm{r}^{(k)} + \bm{s}^{(k)} =
\bm{r}^{(k)} - J^\dagger(\bm{r}^{(k)}) 
F(\bm{r}^{(k)}),
\end{equation*}
where $J^\dagger(\bm{r})$ is the Moore--Penrose pseudoinverse of
$J(\bm{r})$~\cite{bjo96}.

A straightforward computation leads to the expression of the components $f_t$
of $F$
\begin{equation*}
\begin{aligned}
f_t(\bm{r}) = \bm{z}_t^TR^TR\bm{z}_t-1 =& z_{1t}^2r_{11}^2 + 
2z_{1t}z_{2t}r_{11}r_{12} + 2z_{1t}z_{3t}r_{11}r_{13} + 
z_{2t}^2(r_{12}^2+r_{22}^2)\\
& + 2z_{2t}z_{3t}(r_{12}r_{13}+r_{22}r_{23}) + 
z_{3t}^2(r_{13}^2+r_{23}^2+r_{33}^2)-1,
\end{aligned}
\end{equation*}
and, therefore, to the partial derivatives $\partial f_t / \partial 
r_{ij}$. Consequently, the $t$th row of $J$ is
\begin{equation}\label{jacobian}
\begin{aligned}
J_{t,:}(\bm{r}) & =
\begin{bmatrix}
\frac{\partial f_t}{\partial r_{11}} &
\frac{\partial f_t}{\partial r_{12}} &
\frac{\partial f_t}{\partial r_{13}} &
\frac{\partial f_t}{\partial r_{22}} &
\frac{\partial f_t}{\partial r_{23}} &
\frac{\partial f_t}{\partial r_{33}}
\end{bmatrix} \\
& = 2
\begin{bmatrix}
\begin{bmatrix}
r_{11} & r_{12} & r_{13}
\end{bmatrix}\bm{z}_t\bm{z}_t^T & \
\begin{bmatrix}
r_{22} & r_{23}
\end{bmatrix}
\begin{bmatrix} z_{2t}\\ z_{3t} \end{bmatrix}
\begin{bmatrix} z_{2t}& z_{3t} \end{bmatrix}
 & \
r_{33}z_{3t}^2
\end{bmatrix}.
\end{aligned}
\end{equation}

To solve the nonlinear problem \eqref{nonlin} we used the MATLAB function
\texttt{mngn2} discussed in \cite{pr20}, which implements a relaxed version of
the Gauss--Newton method and which is available on the web page
\url{https://bugs.unica.it/cana/software/}.
The algorithm is especially suited for underdetermined nonlinear least-squares
problems, but works nicely also for overdetermined ones.

Formulating the solution of the PS problem with unknown lighting by a nonlinear
model may appear as impractical, when there exists a linear formulation, e.g.,
the Hayakawa procedure.
It has a principal advantage: the positive definiteness of the matrix $G=R^TR$
is not an essential assumption for the application of the method. On the
contrary, if the algorithm is successful in determining $R$, then the matrix
$G$ is positive definite.
This might help in introducing an ideality test for the dataset under scrutiny,
as it will be shown in the following.
Moreover, the computational complexity does not grow excessively, as there are
simple closed formulae for the Jacobian and, in our preliminary experiments,
the Gauss--Newton method proved to converge in a small number of iterations,
when successful.
The performance of the two methods will be compared in the 
numerical experiments described in Section~\ref{sec:expresults} and
discussed in Section~\ref{sec:discussion}.

\section{Dealing with nonideal data}\label{sec:nonideal}

When the light estimation techniques discussed in Sections~\ref{sec:hayakawa}
and~\ref{sec:nonlin} are applied to experimental datasets, usually
computational problems emerge. This is due to the fact that in real
applications some of the assumptions required by the model may be unmet.

The most limiting assumption in the application of the Hayakawa procedure is
the positive definiteness of matrix $G$ in \eqref{systemG}.
For some datasets, the smallest eigenvalue of $G$ turns out to be nonpositive,
leading to a breakdown in the Cholesky algorithm.
This is clearly due to ``nonideality'' of data.


In applicative scenarios, it would be useful to predict if the available data
satisfy the model assumptions by a numerical indicator of ideality.
What is needed is a clear and effective strategy to check this at a reasonable
computational cost.

Our first attempt was to control the error produced by approximating the data
matrix $M$ by the matrix $M_3$, obtained by truncating the factorization
\eqref{factorization} to three terms. In this case, it is known that
$\|M-M_3\|=\sigma_4$, the fourth singular value of $M$, and that such value is
minimal among all matrices of rank 3~\cite{bjo96}.
We tried to vary the composition of the dataset, removing some images from it,
aiming to minimize the value of the approximation error $\sigma_4$.
We also investigated the minimization of the ratio $\rho=\sigma_4/\sigma_3$,
which represents the distance of $M$ from being a rank-3 matrix.
None of these attempts was successful, in the sense that the values of both
$\sigma_4$ and $\rho$ appear to be unrelated to the positive definiteness of
$G$ and to a good 3D reconstruction of the observed object.

Then, we focused on the matrix $H$ defined in \eqref{matrixH}. This matrix
determines the entries of the matrix $G$ via the solution of the system
$H\bm{g}=\bm{1}$, and so is directly related to its spectral properties.
From \eqref{matrixH}, it is immediate to observe that each row of $H$ is
obtained from one column of the matrix $Z$, which is the first tentative
approximation of the light matrix $L$.
So, each row of $H$ is associated with a light direction, and consequently
it depends upon a single image in the dataset. 
Since the matrix $H$ must have at least six rows to ensure a unique solution
matrix $G$, a redundant dataset allows one to investigate the effect of
different image subsets by simply removing from $H$ the rows corresponding to
the neglected images.
One important feature of this process is that the SVD factorization
\eqref{factorization} is performed only once, and the solution of system
\eqref{matrixH} only requires the QR factorization of the relatively small
matrix $H\in\R^{q\times 6}$. 

We propose to validate a dataset using the value of the smallest eigenvalue
$\lambda_3(G)$ of the matrix $G$ as a measure of data ideality.
Indeed, a negative value of $\lambda_3(G)$ is a clear indicator of an
inappropriate shooting technique, and suggests that some images should
be removed from the dataset.
At the same time, to select one between two different reduced datasets, which
are both admissible for the application of the Hayakawa procedure, one may
choose the one which produces the largest eigenvalue $\lambda_3(G)$.
In fact, in our experiments we observed that a small value of $\lambda_3(G)$
may introduce distortions in the 3D reconstruction.

\begin{algorithm}
\caption{Removal of ``unideal'' images from a PS dataset: linear approach.}
\label{algo1}
\begin{algorithmic}[1]
\REQUIRE PS data matrix $M$ of size $p\times q$ ($q$ pictures, each with $p$
pixels)
\ENSURE Set $\cS$ containing the indices of the images to be kept in the
dataset
\STATE Compute the compact SVD $M=U\Sigma V^T$, with
	$V=[\bm{v}_1,\ldots,\bm{v}_q]$, $\bm{v}_i\in\R^q$
\STATE $\widetilde{L}=[\bm{v}_1,\bm{v}_2,\bm{v}_3]^T$
\STATE $\cS=\{1,\ldots,q\}$
\STATE $k=0$, $\mu_0=0$, $L=\widetilde{L}$
\REPEAT
	\STATE $k=k+1$
	\FOR {$i=1,2,3$}\label{Hstart}
		\STATE $H_{:,i}=(L_{i,:}^T)^2$ (componentwise)
	\ENDFOR
	\STATE $H_{:,4}=2(L_{1,:}^T)*(L_{2,:}^T)$ (componentwise)
	\STATE $H_{:,5}=2(L_{1,:}^T)*(L_{3,:}^T)$ (componentwise)
	\STATE $H_{:,6}=2(L_{2,:}^T)*(L_{3,:}^T)$ (componentwise)\label{Hend}
	\FOR {$i=1,\ldots,q$}\label{FORs}
		\STATE Let $\widetilde{H}$ be $H$ with column $i$ removed
		\STATE $\bm{g}=\arg\min\|\widetilde{H}\bm{g}-\bm{1}\|$
		\STATE Construct matrix $G$ from $\bm{g}$
		\STATE $\lambda_i=\min\{\text{eigenvalues of }G\}$
	\ENDFOR\label{FORe}
	\STATE $\lambda_\ell=\max\{\lambda_1,\ldots,\lambda_q\}$
	\IF {$k=1$ and $\lambda_\ell\leq 0$}
		\STATE Stop iteration: unrecoverable breakdown
	\ENDIF
	\STATE Store in $t$ the $\ell$th element of $\cS$
	\STATE Remove $t$ from $\cS$\label{Srem}
	\STATE $q=q-1$
	\STATE $\mu_k=\lambda_\ell$
	\IF {fast version}
		\STATE $L=\widetilde{L}_{:,\cS}$ ($L$ contains the columns of
			$\widetilde{L}$ with indices in $\cS$)\label{econ}
	\ELSE 
		\STATE Remove column $\ell$ from $M$\label{updstart}
		\STATE Compute the compact SVD $M=U\Sigma V^T$, with
			$V=[\bm{v}_1,\ldots,\bm{v}_q]$
		\STATE $L=[\bm{v}_1,\bm{v}_2,\bm{v}_3]^T$\label{updstop}
	\ENDIF
\UNTIL {$\mu_k<\mu_{k-1}$ or $q=6$}\label{stop1}
\STATE $\cS=\cS\cup\{t\}$
\end{algorithmic}
\end{algorithm}

This idea inspired the numerical method outlined in Algorithm~\ref{algo1}.
It stands on the assumption that in the given dataset of $q$ images there
exists a subset of $q-1$ images which leads to a positive definite matrix $G$.
If this condition is not met, then some preprocessing is needed to exclude the
worst images from the dataset.

At the beginning, we construct the rank-3 factor $L_3$ from
\eqref{factorization}, here denoted as $\widetilde{L}$, and initialize the
indices set $\cS$ to $\{1,\ldots,q\}$.
Inside the main loop, at lines \ref{Hstart}--\ref{Hend}, the matrix $H$ with
rows indexed in $\cS$ is constructed.
Then, each row is iteratively removed from $H$ and the smallest eigenvalue of
the corresponding matrix $G$ is stored; see lines \ref{FORs}--\ref{FORe}.
The largest among the found eigenvalues identifies the index to be removed from
$\cS$ (line \ref{Srem}), that is, the image to be excluded.
At this point, the corresponding column is removed from $M$ and its compact SVD
decomposition is recomputed; see lines \ref{updstart}--\ref{updstop}.
The iteration continues until the sequence of selected eigenvalues stops
increasing, or the number of images $q$ falls below six.

We explored the possibility of reducing the computational cost by avoiding
recomputing the SVD decomposition of $M$, and extracting the matrix $L$ at
each iteration from the initial factorization; see line~\ref{econ}. The
resulting algorithm, denoted as ``fast'' in Algorithm~\ref{algo1}, is
less accurate, as it will be shown in the following.

\begin{algorithm}
\caption{Removal of ``unideal'' images from a PS dataset: nonlinear approach.}
\label{algo2}
\begin{algorithmic}[1]
\REQUIRE PS data matrix $M$ of size $p\times q$ ($q$ pictures, each with $p$
pixels)
\ENSURE Set $\cS$ containing the indices of the images to be kept in the
dataset
\STATE $\cS=\{1,\ldots,q\}$
\STATE $k=0$, $\rho_0=0$
\REPEAT
	\STATE $k=k+1$
	\IF {fast version}
		\STATE Let $\widetilde{M}$ contain the columns of $M$ indexed
			in $\cS$\label{esvd2s}
		\STATE Compute the compact SVD $\widetilde{M}=U\Sigma V^T$,
			with $V=[\bm{v}_1,\ldots,\bm{v}_q]$
		\STATE $\widetilde{L}=[\bm{v}_1,\bm{v}_2,\bm{v}_3]^T$
			\label{esvd2e}
	\ENDIF
	\FOR {$i=1,\ldots,q$}
		\IF {fast version}
			\STATE Let $L$ be $\widetilde{L}$ with column $i$
				removed\label{unupL}
		\ELSE
			\STATE Let $\widetilde{M}$ be $M$ with column $i$
				removed\label{start2}
			\STATE Compute the compact SVD 
				$\widetilde{M}=U\Sigma V^T$,
				with $V=[\bm{v}_1,\ldots,\bm{v}_q]$
			\STATE $L=[\bm{v}_1,\bm{v}_2,\bm{v}_3]^T$\label{end2}
		\ENDIF
		\STATE Run \texttt{mngn2}~\cite{pr20} on problem \eqref{nonlin}
			and obtain $\eta_i=\gamma_6/\gamma_5$ at
			convergence\label{mngn}
	\ENDFOR
	\STATE $\eta_\ell=\max\{\eta_1,\ldots,\eta_q\}$
	\STATE Store in $t$ the $\ell$th element of $\cS$
	\STATE Remove $t$ from $\cS$\label{cS2}
	\STATE $q=q-1$
	\STATE $\rho_k=\eta_\ell$
\UNTIL {$\rho_k<\rho_{k-1}$ or $q=6$}\label{stop2}
\STATE $\cS=\cS\cup\{t\}$
\end{algorithmic}
\end{algorithm}

The nonlinear model for identifying the lights location, outlined in
Section~\ref{sec:nonlin}, also leads to a procedure for excluding unideal
images from a redundant dataset.
Indeed, we noticed that in the presence of images that deviate from ideality
and that would lead to a nonpositive definite matrix $G$ in the Hayakawa 
procedure, the Gauss--Newton method diverges, producing a meaningless solution.
The reason is that the Jacobian \eqref{jacobian} becomes rank deficient as the
iterates reach a neighborhood of the solution.
Indeed, the $q\times 6$ matrix $J(\bm{r}^{(k)})$ ($q\geq 6$) should have rank
6, but as the iteration progresses it may happen that its numerical rank falls
below this value.
Denoting by $\gamma_i$, $i=1,\ldots,6$, the singular values of
$J(\bm{r}^{(k)})$, we propose to use the ratio between the sixth and the fifth
singular values 
\begin{equation}\label{ratio}
\eta=\frac{\gamma_6}{\gamma_5}
\end{equation}
as an indicator of nonideality.
This leads to the numerical procedure reported in Algorithm~\ref{algo2}.

After initializing the set $\cS$ to the column indices of the data matrix $M$,
the main loop starts. There, the algorithm removes one column from $M$,
computes its SVD decomposition, and constructs a tentative matrix $L$; see
lines \ref{start2}--\ref{end2}.
Then (line~\ref{mngn}), the nonlinear least-squares method \texttt{mngn2}
from~\cite{pr20} is applied to problem \eqref{nonlin}.
The function returns the value of the ratio $\eta$ in \eqref{ratio} for the
Jacobian evaluated at the converged iteration.
When all the columns have been analyzed, the index of the maximum value of
$\eta$ identifies the best column configuration, so it is removed from the set
$\cS$ (line \ref{cS2}) and the iteration is restarted.
The computation ends when the sequence of ratios stops increasing, or when just
six images are left in the dataset.

Algorithm~\ref{algo2} also contains a ``fast'' version of the method.
To reduce the complexity, the SVD factorization is not computed at each
step of the \textsc{for} loop, but just before the loop starts (lines
\ref{esvd2s}--\ref{esvd2e}).
Then, the matrix $L$ is constructed by selecting the relevant
columns from the unupdated factor $V$ (line \ref{unupL}).
This version of the method produces slightly different results. Its performance
will be discussed in the numerical experiments.

We remark that the computational cost of Algorithms~\ref{algo1} and~\ref{algo2}
does not usually have an impact on real-time processing.
Indeed, in most applications a dataset is analyzed only once, right after
acquisition, and a reduced dataset is then stored for further processing.

\section{Numerical experiments}\label{sec:expresults}

To investigate the effectiveness of Algorithms~\ref{algo1} and~\ref{algo2}, we
implemented them in the MATLAB programming language. The two functions are
available from the authors upon request and will be included in the MATLAB
toolbox presented in~\cite{cdfrv21}. The same toolbox, available at the web
page \url{https://bugs.unica.it/cana/software/}, was also used to compute the
3D reconstructions.

\begin{figure}[htb]
\centerline{\includegraphics[width=.5\linewidth]{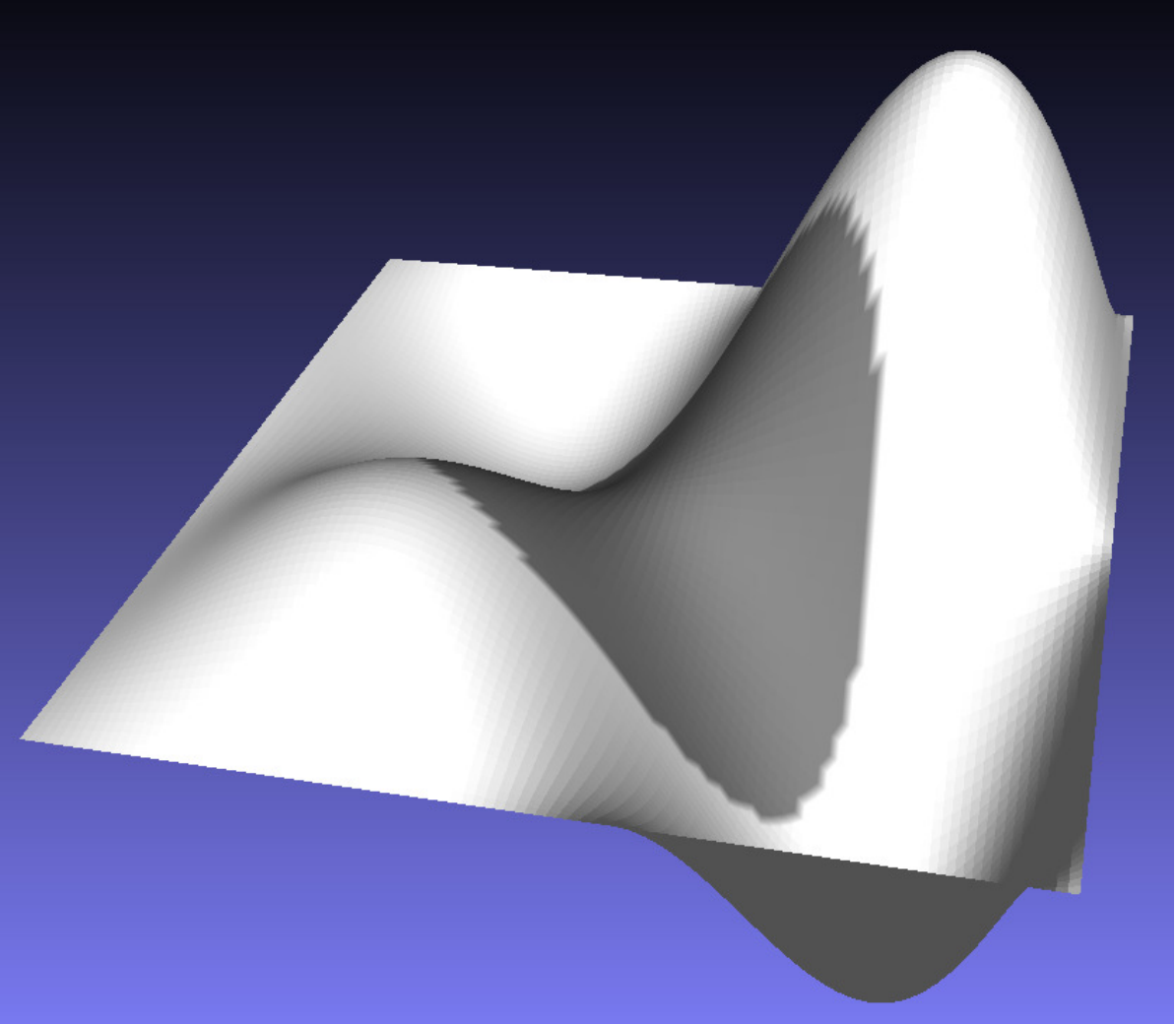}\phantom{XX}
\includegraphics[width=.43\linewidth]{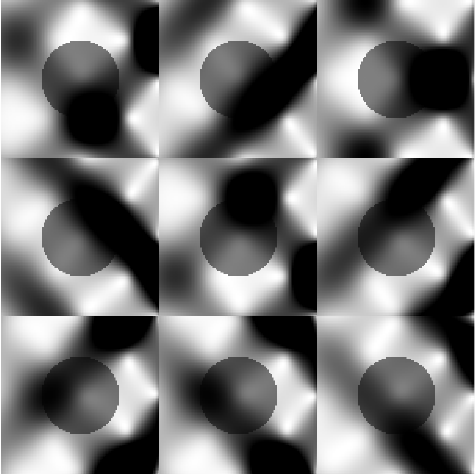}}
\caption{The synthetic dataset used in the experiments is generated by applying
the PS forward model to the surface displayed on the left, with a
discretization of size $101\times 101$; the images reported on the right are
obtained by choosing 9 different light orientations.}
\label{fig:syntimages}
\end{figure}

We initially apply the two methods to a synthetic dataset, first considered
in~\cite{cdfrv21}. The analytical expression $z=u(x,y)$ of the observed
surface, displayed in the left-hand pane of Figure~\ref{fig:syntimages}, is
fixed a priori, and a simple white and gray image is texture mapped to it.
The synthetic images are obtained by applying the PS forward model to the
surface. 
The images displayed on the right of Figure~\ref{fig:syntimages} were
produced by choosing a set of nine lighting directions.
In this case, the light sources were located at an infinite distance from the
target, but the model allows choosing the lighting distance using the width
$A$ of the domain $\Omega$ (see \eqref{discretization}) as a unit of measure.

We used the dataset reported in Figure~\ref{fig:syntimages}, with lights at
infinity and no noise, substituting the third image with one obtained with
a light source at distance $\delta A$, for $\delta=10,8,6,4,2$. We also
contaminated this image with Gaussian noise having mean value zero and standard
deviation $0.1$.
This produces a dataset where the only unideal image is the number 3.

\begin{figure}[htb]
\centerline{\includegraphics[width=.7\linewidth]{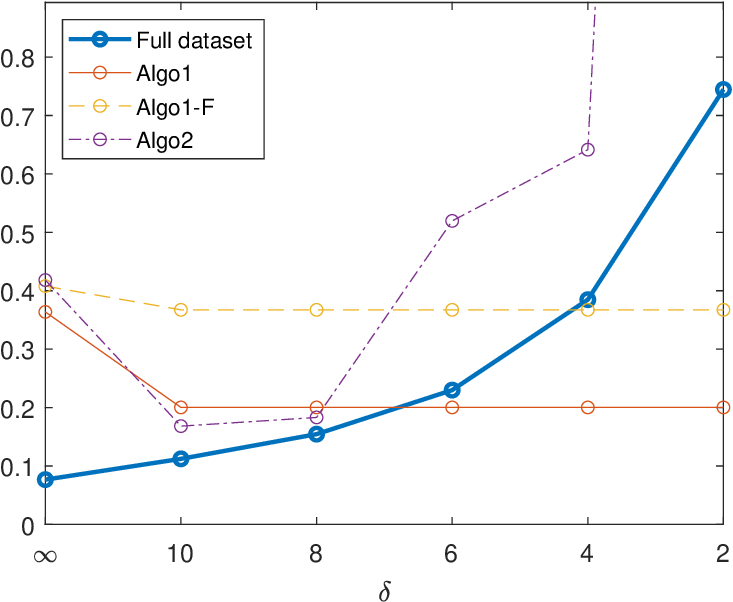}}
\caption{Relative errors in the $\infty$-norm, obtained by processing either
the whole dataset or the reduced datasets suggested by the algorithms under
scrutiny; $\delta A$ is the distance of the light source from the object.}
\label{fig:syntherror}
\end{figure}

The blue thick line in Figure~\ref{fig:syntherror} represents the relative
error
$$
E(\delta) = \frac{\|u(x,y)-\tilde{u}(x,y)\|_\infty}{\|u(x,y)\|_\infty},
$$ where 
$\|u(x,y)\|_\infty=\max_{(x,y)\in\Omega}|u(x,y)|$,
$u(x,y)$ is the model surface, and $\tilde{u}(x,y)$ is the
reconstruction obtained by processing the whole dataset, where the third
light source is at distance $\delta A$.
The other three lines represent the error obtained by processing a reduced
dataset, according to the recommendation of Algorithm~\ref{algo1} (Algo1), the
fast version of the same method (Algo1-F), and Algorithm~\ref{algo2}
(Algo2).

It can be seen that, while the error corresponding to the whole dataset
increases as $\delta$ decreases, the two versions of Algorithm~\ref{algo1}
identify a subset for which the error is reduced, when the light source is very
close. Indeed, both versions select image 3 as the first to be removed,
while the fast version (incorrectly) also excludes image 1 from the dataset.
On the contrary, Algorithm~\ref{algo2} fails, as it excludes images 5 and 6,
leading to a larger error. The fast version of Algorithm~\ref{algo2} leads to
even worse results.

\begin{figure}[htb]
\centerline{\includegraphics[width=.8\linewidth]{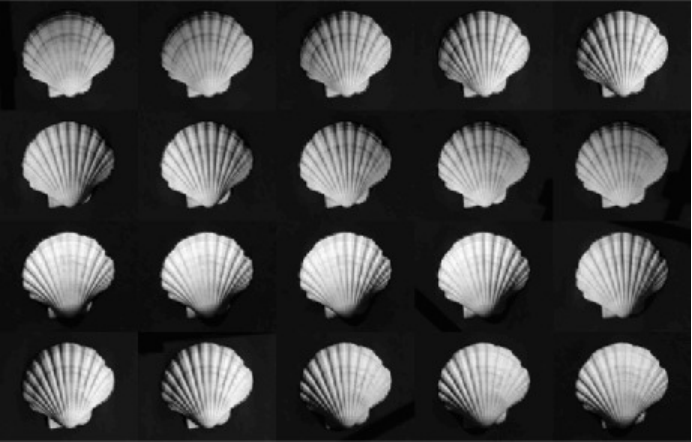}}
\caption{The \emph{Shell3} dataset is composed of 20 pictures of a seashell
illuminated by sunlight; the light direction rotates counterclockwise;
see~\cite{rrmyb23} for details.} \label{fig:shell}
\end{figure}

While the first experiment was concerned with the case of a perfect Lambertian
reflector illuminated by a close light source, in the second one we investigate
the case of a non-Lambertian surface with lights at almost infinite distance
from it.
The \emph{Shell3} dataset, displayed in Figure~\ref{fig:shell}, was used for
the first time in~\cite{rrmyb23}.
A preliminary version of the dataset has been considered in~\cite{cdfrv21}.
The pictures were obtained by placing a seashell, approximately 10 cm wide, on a
rotating platform in the open air. The same platform holds a tripod with a
camera observing the shell from above, at a distance of about 1 m;
see~\cite{rrmyb23}.
Twenty images were taken by letting the platform rotate clockwise under direct
sunlight.

The relatively large distance between the camera and the shell, coupled with a
lens with a focal length of 85 mm, produces a reasonable approximation of an
orthographic projection.
The light source is virtually at an infinite distance from the object.
Nevertheless, the Hayakawa procedure fails when it is applied to the whole
dataset, as in this case the matrix $G$ in \eqref{systemG} is not positive
definite.
This is probably due to the fact that the shell surface does not reflect the
light according to Lambert's law.

We report in the following the indices of the images that were excluded from
the dataset, after processing the data matrix $M$ of size $702927\times 20$, by
each of the 4 methods discussed in Section~\ref{sec:nonideal}:
\begin{description}
\item[Algorithm~\ref{algo1}:] 3, 13, 2;
\item[Algorithm~\ref{algo1}-fast version:] 3, 13, 2, 17, 20, 5, 11;
\item[Algorithm~\ref{algo2}:] 3, 20, 7, 2;
\item[Algorithm~\ref{algo2}-fast version:] 3, 20, 5.
\end{description}
{Each of the four lines contains the ordered list of the images identified
at each iteration of Algorithms~\ref{algo1} (line~\ref{Srem})
and~\ref{algo2} (line~\ref{cS2}). Each list is terminated when the stop
condition of the main loop is met (lines~\ref{stop1} and~\ref{stop2},
respectively).}

\begin{figure}[htb]
\includegraphics[width=.49\linewidth]{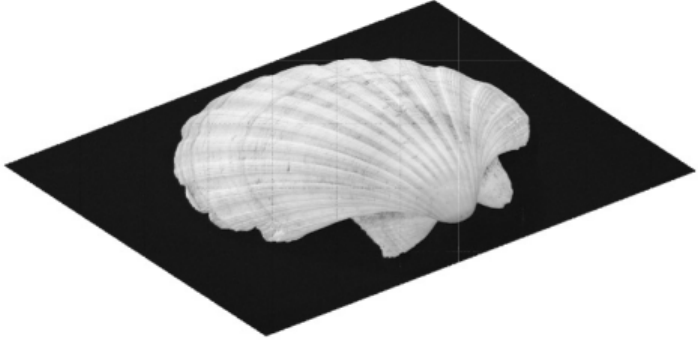}
\includegraphics[width=.49\linewidth]{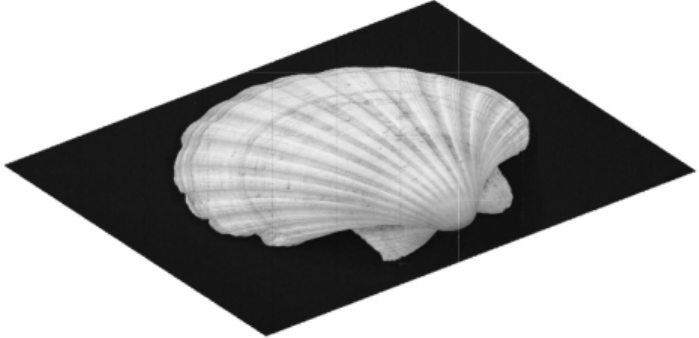}
\includegraphics[width=.49\linewidth]{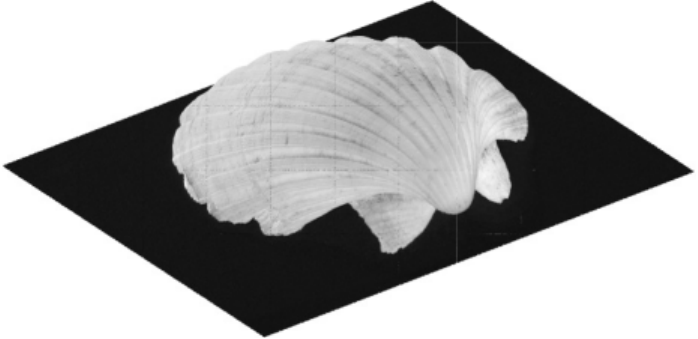}
\includegraphics[width=.49\linewidth]{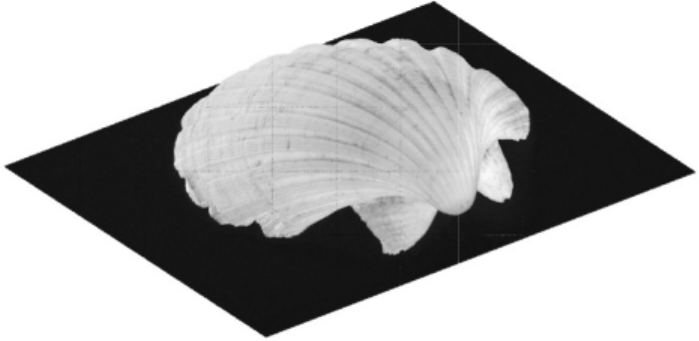}
\caption{3D reconstructions obtained after removing from the \emph{Shell3}
dataset the pictures selected by Algorithm~\ref{algo1} (top-left), the fast
version of Algorithm~\ref{algo1} (top-right), Algorithm~\ref{algo2}
(bottom-left), and the fast version of Algorithm~\ref{algo2} (bottom-right).}
\label{fig:shell-results}
\end{figure}

\begin{figure}[htb]
\centerline{\includegraphics[width=.6\linewidth]{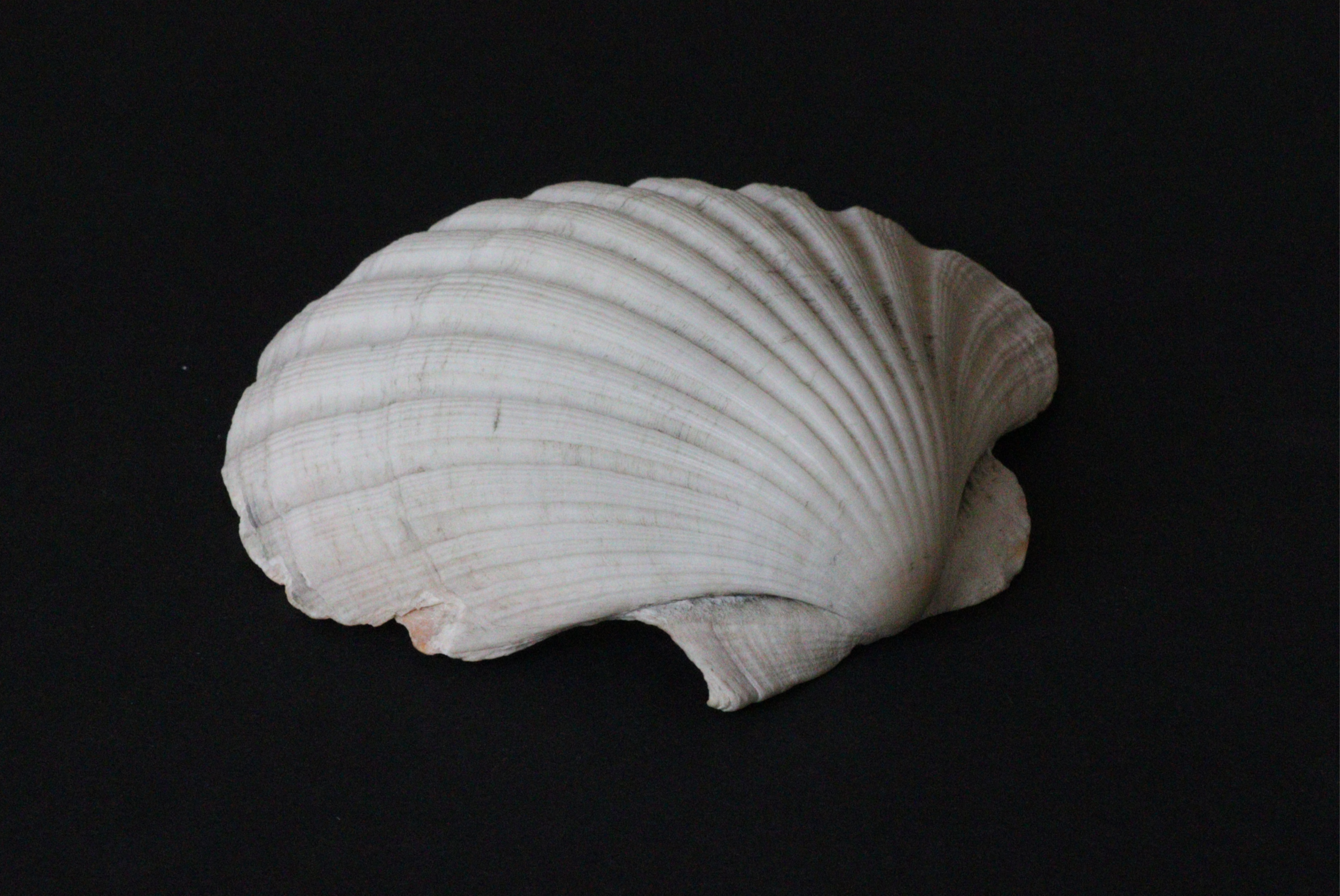}}
\caption{Color picture of the seashell from which the \emph{Shell3} dataset was
created.}
\label{fig:shell-image}
\end{figure}

Figure~\ref{fig:shell-results} reports a view of the reconstructions obtained
by the four methods. They should be compared with the picture of the real
seashell, displayed in Figure~\ref{fig:shell-image}.
Even if in this case it is not possible to obtain a numerical measure of the
error, it is evident that the surfaces recovered by Algorithm~\ref{algo2} are
much closer to the original, while those produced by Algorithm~\ref{algo1} are
considerably flatter.

\section{Discussion}\label{sec:discussion}

Our numerical experiments faced two situations of particular importance in PS.
In the first one, all the images in the synthetic dataset meet the assumptions
of the PS model, except for one image, which represents the same Lambertian
surface lighted by a source approaching the target and is affected by $10\%$
Gaussian noise.

As Figure~\ref{fig:syntherror} shows, while processing the whole dataset leads
to an error which increases as the light source approaches the surface,
Algorithm~\ref{algo1} is able to identify the disturbing image, and its removal
leads to a substantial reduction in the error. The fast version of the same
method proves to be less accurate, but still produces a smaller error than the
unreduced dataset.
Algorithm~\ref{algo2}, on the contrary, is not able to select the unideal
image, and selects a subset that degrades the quality of the reconstruction.

Figure~\ref{fig:syntherror} also shows that for $\delta\geq 8$ the error in
the reconstruction is not much larger than with a light source at infinity, and
removing the single unideal image from the dataset produces slightly worse
results.
The quality of the recovered surface is strongly affected by a close light
source only when $\delta\leq 4$, that is, in very narrow shooting scenarios.

The second dataset is experimental and it is correctly illuminated, since the
sunlight produces practically parallel rays. Anyway, the surface of the
seashell is probably not Lambertian and the pictures are affected by measuring
errors.
The standard Hayakawa procedure is not applicable to the dataset, as it breaks
when computing the Cholesky factorization of $G$ in \eqref{systemG}, so the
dataset must be reduced.
We see that both developed algorithms are able to produce a meaningful
solution, but the reconstruction produced by Algorithm~\ref{algo2} seems to
better approximate the original seashell, even when the fast version of the
method is employed.
{These results are in contrast with the graph in
Figure~\ref{fig:syntherror}, which reports the error behaviour for the
synthetic dataset, where Algorithm~\ref{algo1} performs better.}

{The above consideration suggests that both methods may have a role in
detecting a lack of ideality in a large dataset, for what regards the presence
of close light sources and of a non-Lambertian target. Our further studies will
address the development of a procedure which blends the two approaches for
analyzing a redundant dataset.
An immediate advantage of the two algorithms is that they suggest an ordered
list of pictures to be excluded from the numerical computation, that may assist
the user in preprocessing a dataset to better suit his/her needs.}
This may be accomplished by visually comparing the reconstructed surfaces to the
true object. Indeed, the computing time required by the algorithms
proposed in this paper to reduce the dataset, and by the package presented in
\cite{cdfrv21} to generate the reconstructions, allows for fast on-site
processing on a laptop.

To conclude, the proposed new approaches look promising for determining a
subset of an experimental PS dataset that better approximates the strict
assumptions required by Lambert's model under unknown lighting.
One of their limitations is that they always exclude at least one picture from
the initial dataset, even if it is a perfect one. This means that they are not
able to evaluate a dataset in itself, but only to compare successive subsets of
the initial collection of images, and that their application must always
be supervised by a human operator. We believe that a hybrid method, which
keeps into account the forecasts of both algorithms, may improve their
performance, but this requires more work.
What is also needed is a wide experimentation on well-known collections of PS
data, such as~\cite{shi2016benchmark,ren2022diligent102}, as well as
on images acquired in less restrictive conditions than the ones considered in
this paper, to compare the performance of the proposed methods to other
available techniques.
This will be the subject of future research.

\medskip

\section*{Acknowledgements}
This work was partially supported 
by Fondazione di Sardegna, Progetto biennale bando 2021, ``Computational
Methods and Networks in Civil Engineering (COMANCHE)'',
the INdAM-GNCS 2023 project ``Tecniche numeriche per lo studio dei problemi
inversi e l'analisi delle reti complesse'',
E. C. gratefully acknowledges CeSIM for the financial support provided under
the project ``Documentation techniques for rock and wall engravings found in
Domus de Janas, caves, shelters and open rocks''.

\bibliographystyle{siam}
\bibliography{bibphoto3d}

\end{document}